\documentclass{gtpart}
\usepackage{amscd,amsmath,amssymb}

\title[Immersions of $RP^{2^e-1}$]{Immersions of $RP^{2^e-1}$}

\author{Donald M. Davis}
\givenname{Donald}
\surname{Davis}
\email{dmd1@lehigh.edu}
\address{Lehigh University\\Bethlehem, PA 18015}
\urladdr{http://www.lehigh.edu/~dmd1}

\author{Giora Dula}
\givenname{Giora}
\surname{Dula}
\email{giora@mail.netanya.ac.il}
\address{Netanya College\\Netanya 42365, Israel}
\urladdr{}

\author{Jes\'us Gonz\'alez}
\givenname{Jes\'us}
\surname{Gonz\'alez}
\email{jesus@math.cinvestav.mx}
\address{Departamento de Matem\'aticas\\ Centro de Investigaci\'on y de Estudios Avanzados del IPN\\ M\'exico DF 07000, M\'exico}
\urladdr{}

\author{Mark Mahowald}
\givenname{Mark}
\surname{Mahowald}
\email{mark@math.northwestern.edu}
\address{Northwestern University\\Evanston, IL 60208}
\urladdr{}


\arxivreference{AT/0701327}  
\arxivpassword{f5pxv}   

\newtheorem{thm}[equation]{Theorem}
\numberwithin{equation}{section}

\newtheorem{lem}[equation]{Lemma}

\newtheorem{defin}[equation]{Definition}

\newtheorem{prop}[equation]{Proposition}

\newtheorem{rmk}[equation]{Remark}

\renewcommand{\theequation}{\thesection.\arabic{equation}}

\def\mapright#1{\ \smash{\mathop{\longrightarrow}\limits^{#1}}\ }

\def\ss{\smallskip}

\def\la{\langle}
\def\ra{\rangle}
\def\on{\operatorname}

\def\lbar{{\overline{\ell}}}
\def\qed{\quad\rule{8pt}{8pt}\bigskip}

\def\a{\alpha}

\def\eps{\epsilon}

\def\bx{{\mathbf x}}

\def\dbar{{\overline d}}

\def\tfrac{\textstyle\frac}

\def\w{\wedge}
\def\equ{\begin{equation}}
\def\b{\beta}

\def\endeq{\end{equation}}

\def\coker{\on{coker}}

\def\gt{\widetilde{L}}

\def\gd{{\on{gd}}}

\def\R{{\mathbb R}}

\begin{document}

\dedicatory{This paper is dedicated to Michael Barratt
on the occasion of his 81st birthday.}

\begin{abstract}
We prove that $RP^{2^e-1}$ can be immersed in $\R^{2^{e+1}-e-7}$ provided $e\ge 7$. If $e \ge 14$, this is 1 better than previously known
immersions. Our method is primarily an induction on geometric dimension, with compatibility of liftings being a central issue.
\end{abstract}

\maketitle

\section{Statement of result and background}\label{intro}

Our main result is the following immersion theorem for real projective
spaces.

\begin{thm}
If $e\ge7$, then $RP^{2^e-1}$ can be immersed in $\R^{2^{e+1}-e-7}$. \label{mainthm}\end{thm}

This improves, in these cases, by 1 dimension upon the result of Milgram
\cite{Mil}, who proved, by constructing bilinear maps, that if $n\equiv7$
mod 8,
then $RP^n$ can be immersed in $\R^{2n-\a(n)-4}$, where $\a(n)$ denotes the
number of 1s in the binary
expansion of $n$. In \cite[Theorem 1.2]{DM}, the first and fourth authors used
obstruction theory to prove that if $n\equiv7$ mod 8, then
$RP^n$ can be immersed in $\R^{2n-D}$, where $D=14,16,17,18$ if
$\a(n)=7,8,9,\ge10$. That result, with $n=2^e-1$, is
stronger than ours for $ e\le 12$. If $e\ge13$, then our result
improves on the result of \cite{DM}
by $e-13$ dimensions.
Thus \fullref{mainthm} improves on all known results by 1 dimension if
$e\ge14$.

In \cite{J}, James proved that $RP^{2^e-1}$ cannot be immersed in
$\R^{2^{e+1}-2e-\delta}$ where $\delta=3,2,2,4$ for
$e\equiv 0,1,2,3$ mod 4. In \cite{GM}, Gitler and  Mahowald
announced an immersion result for $RP^{2^e-1}$
in dimension 1 greater than
that of James' nonimmersion, which would have been optimal. However, a
mistake in the argument of \cite{GM} was
pointed out by Crabb and Steer. The approach of our paper was initiated
by Mahowald around 1970 in an unpublished attempt to prove an optimal immersion of
$RP^{2^e-1}$. In order to improve our result to this extent, we would need
to show compatibility of our liftings with liftings given by the
Radon-Hurwitz theorem (\cite{Eck}).

\section{Outline of proof}\label{outlinesec}
In this section we outline the proof of \fullref{mainthm}. In subsequent
sections, we will fill in details.

If $\theta$ is a vector bundle over a compact connected space $X$, we
define the geometric dimension of $\theta$, denoted $\gd(\theta)$, to be
the fiber dimension of $\theta$ minus the maximum number of linearly
independent sections of $\theta$. Equivalently, if $\dim(\theta)=n$, then
$\gd(\theta)$ equals the smallest integer $k$ such that the map
$X\mapright{\theta} BO(n)$ which classifies $\theta$ factors through
$BO(k)$.
The following lemma is standard (See eg Sanderson \cite[Theorem 4.2]{San}).
Here and throughout, $\xi_n$ denotes the Hopf line bundle over $RP^n$. We
will often write $P^n$ instead of $RP^n$, and will denote
the stunted space $P^n/P^{k-1}$ as $P^n_k$.
\begin{lem} Let $\phi(n)$ denote the number of positive integers $i$
satisfying $i\le n$ and $i\equiv0,1,2,4 \mod 8$.
Suppose $n>8$. Then $RP^n$ can be immersed in $\R^{n+k}$ if and only if
$\gd((2^{\phi(n)}-n-1)\xi_n)\le k$.\end{lem}
Thus \fullref{mainthm} will follow from the following result, to the proof
of which the remainder of this
paper will be devoted.

\begin{thm}\label{gdthm} If $e\ge7$, then
$\gd((2^{2^{e-1}-1}-2^e)\xi_{2^e-1})\le 2^e-e-6$.\end{thm}

The bulk of the work toward proving \fullref{gdthm} will be a determination
of upper bounds
for $\gd(2^e\xi_n)$ for all $n\equiv7$ mod 8 by induction on $e$, starting with $e=7$. A similar
method could be employed for
all $n$, but we restrict to $n\equiv7$ mod 8 to simplify the already
formidable arithmetic. We let $A_k=RP^{8k+7}$,
and denote $\gd(m\xi_{8k+7})$ by $\gd(m,k)$.

The classifying map for $2^e\xi_{8k+7}$ will be
viewed as the following composite.
\begin{equation}\label{comp}A_k\mapright{d} (A_k\times
A_k)^{(8k+7)}\hookrightarrow\bigcup_j A_j\times A_{k-j}\mapright{f\times f}
BO_{2^{e-1}}\times BO_{2^{e-1}}\to BO_{2^e}.\end{equation}
Here $d$ is a cellular map homotopic to the diagonal map, $X^{(n)}$ denotes
the $n$--skeleton of $X$, and $f$ classifies
$2^{e-1}\xi$.
We write $BO_m$ for $BO(m)$ for later notational convenience.

As a first step, we would like to use \eqref{comp} to deduce that
$$\gd(2^e,k)\le\max\{\gd(2^{e-1},j)+\gd(2^{e-1},k-j):0\le j\le k\}.$$
In order to make this deduction, we need to know that the liftings of the
various $2^{e-1}\xi_{8j+7}$ to various
$BO_m$ have been made compatibly.

\begin{defin} If $\theta$ is a vector bundle over a filtered space
$X_0\subset\cdots\subset X_k$,
we say that $$\gd(\theta|X_{i})\le d_i\text{ compatibly for }i\le k$$ if
there is a commutative diagram
$$\begin{CD} X_{0}@>>> X_{1} @>>>\cdots@>>>X_{k}@.{}\\
@VVV @VVV @. @VVV@.\\
BO_{d_0}@>>> BO_{d_1} @>>> \cdots @>>> BO_{d_k}@>>>
BO_{\dim(\theta)}\label{compatdef}\end{CD}$$
where the map $X_k\to BO_{\dim(\theta)}$ classifies $\theta$, and the
horizontal maps are the usual inclusions.\end{defin}

\begin{rmk}
{\rm In our filtered spaces, we always assume that the inclusions are cofibrations.}
\end{rmk}

\begin{rmk}
{\rm Isomorphism classes of $n$--dimensional vector bundles over $X$ correspond to homotopy classes of maps of $X$ into $BO_n$. Thus one would
initially say that the diagram in \fullref{compatdef} commutes up to homotopy. However, by \fullref{xx}, we may interpret this diagram, and
other homotopy commutative diagrams that occur later, as being strictly commutative. To apply the lemma, we will often, at the outset, replace
maps $BO_n\to BO_{n+k}$ by homotopy equivalent fibrations.}
\end{rmk}

\begin{lem}\label{xx} If
$$\begin{CD} A@>f>> E\\ @ViVV @VVpV\\
X@>>g> B\end{CD}$$ commutes up to homotopy and $p$ is a fibration, then
$f$ is homotopic to a map $f'$ such that $p\circ f'=g\circ i$.\end{lem}

\begin{proof} Let $H \co A\times I\to B$ be a homotopy from $p\circ f$ to
$g\circ i$.
By the definition of fibration, there exists $\widetilde H \co A\times I\to E$
such that $p\circ\widetilde H=H$ and $\widetilde H|A\times 0=f$. Then
$\widetilde H|A\times 1$ is our desired $f'$.\end{proof}

If $X_0\subset\cdots\subset X_k$
and $Y_0\subset\cdots\subset Y_k$ are filtered spaces, we define, for $0\le
i\le k$,
$$(X\times Y)_i:=\bigcup_{j=0}^i X_j\times Y_{i-j}.$$
Then $(X\times Y)_0\subset\cdots\subset(X\times Y)_k$ is clearly a filtered space. We will prove the following general result in Section
\fullref{compatsec}.

\begin{prop}\label{compatprop} Suppose $\gd(\theta|X_{i})\le d_i$
compatibly for $i\le k$ and
$\gd(\eta|Y_{i})\le d'_i$ compatibly for $i\le k$. For $0\le j\le k$, let
$e_j=\max(d_i+d'_{j-i}: 0\le i\le j)$. Then
$\gd(\theta\times\eta|(X\times Y)_{j})\le e_j$ compatibly for $j\le k$.
Moreover, if $X=Y$ and $\theta=\eta$, then the maps $(X\times
X)_{j}\mapright{f} BO_{e_j}$ can be chosen to satisfy $f\circ T=f$, where
$T \co X\times X\to X\times X$ interchanges factors.\end{prop}

We will begin an induction by deriving in \fullref{128prop} some compatible bounds for $\gd(128,i)$. Proposition \fullref{compatprop} will,
after restriction under the diagonal map, allow us to prove $\gd((\sum 2^{e_i})\xi_n)\le\max\{\sum\gd(2^{e_i}\xi_{m_i}):\sum m_i=n\}$. These
bounds are not yet strong enough to yield new immersion results.
We must improve the bounds by taking advantage of paired obstructions. The following result will be proved in Section \fullref{compatsec}.

\begin{prop}\label{pairedprop} Let $BO_n[\rho]$ denote the pullback of
$BO_n$ and the $(\rho-1)$-connected cover $BO[\rho]$ over $BO$, and let
$s=\min(\rho+2m-1,4m-1)$.
\begin{enumerate}
\item There are equivalences $c'_1$ and $c'_2$ such that the following
diagram commutes.
$$\begin{CD} BO_{2m}[\rho]^{(s)}@>q_1>>
(BO_{2m}[\rho]/BO_{2m-1}[\rho])^{(s)}@>c'_1>>S^{2m}\\
@Vp_2VV @Vp_2'VV @ViVV\\
BO_{2m+1}[\rho]^{(s)}@>q_2>>
(BO_{2m+1}[\rho]/BO_{2m-1}[\rho])^{(s)}@>c'_2>>\Sigma P^{2m}_{2m-1}.
\end{CD}$$
Preparatory to the next two parts, we expand this diagram as follows, with
$c_i=c_i'\circ q_i$ and $(X,A)$ a finite
CW pair.
$$\begin{CD} A@>f_1>>BO_{2m-1}[\rho]^{(s)}@.{}\\
@VjVV @Vp_1VV @.\\
X@.BO_{2m}[\rho]^{(s)}@>c_1>> S^{2m}\\
@. @Vp_2VV @ViVV\\
{}@. BO_{2m+1}[\rho]^{(s)} @>c_2>>\Sigma P_{2m-1}^{2m}.
\end{CD}$$
\item Suppose $\dim(X)<s$, and we are given $X\mapright{f}
BO_{2m}[\rho]^{(s)}$ such that $f\circ j=p_1\circ f_1$ and $c_1\circ f$
factors as
$X\to X/A\mapright{g} S^{2m}$ with $[g]$ divisible by $2$ in
$[X/A,S^{2m}]$.\footnote{Note that $[X/A,S^{2m}]$ is in the stable range,
from
which it gets its group structure.} Then $p_2\circ f$ lifts to a map
$X\mapright{\ell}BO_{2m-1}[\rho]^{(s)}$ whose restriction to $A$ equals
$f_1$.
\item Suppose, on the other hand, $\dim(X)\le s$, and we are given
$X\mapright{f'} BO_{2m+1}[\rho]^{(s)}$ such that $f'\circ j=p_2\circ
p_1\circ f_1$ and $c_2\circ f'$ factors as $X\to X/A\mapright{g'}\Sigma
P_{2m-1}^{2m}$ with $[\Sigma g']$ divisible by $2$ in the stable group
$[\Sigma X/A,\Sigma^2 P_{2m-1}^{2m}]$. Then $f'$ is homotopic rel $A$ to a
map which lifts to $BO_{2m}[\rho]^{(s)}$.
\end{enumerate}
\end{prop}

In Section \fullref{indnsec}, we will implement \fullref{compatprop} and \fullref{pairedprop} to prove that the last part of
the following important result follows by induction on $e$ from the first five parts and its validity when $e=7$, while in Section
\fullref{numbersec}, we will establish the first five parts.

\begin{thm}\label{gthm} There is a function $g(e,k)$ defined for $e\ge7$
and $k\ge0$ satisfying the following.
\begin{enumerate}
\item If $k\ge2^{e-3}$, then $g(e,k)=2^e$.
\item For all $e$, $g(e,0)=g(e,1)=0$, and, if $2\le k\le 2^e$, then  $g(e,k)\ge 4k+4$.
\item If $0\le \ell\le k/2$, then $g(e+1,k)\ge g(e,\ell)+g(e,k-\ell)-1$.
\item If, for some $\ell$ with $0\le\ell\le k/2$, we have $g(e+1,k)=g(e,\ell)+g(e,k-\ell)-1$, then, for all $\ell$ with $0\le\ell\le (k-1)/2$, we have
$g(e,\ell)+g(e,k-1-\ell)<g(e+1,k)$ and, if also $k$ is even, then $g(e+1,k)\ge 2g(e,k/2)+1$.
\item For all $e$ and $k$, $g(e,k)\ge g(e,k-1)$. \item $\gd(2^e,k)\le g(e,k)$ compatibly for all $k$.
\end{enumerate}\end{thm}

The function $g$ will be defined in  \eqref{complement} and \fullref{geq}. In \fullref{gtbl}, we list its values for
small values of the parameters. We prefer not to tabulate the values $g(e,k)=2^e$ when $k>2^{e-3}$.

\begin{table} [ht!]
\begin{center}
\small{
\begin{tabular}{lr|rrrrrrrrrrrrrrrr}
&&&&&&$k$&&&&&&&&&&&\\
&&$1$&$2$&$3$&$4$&$5$&$6$&$7$&$8$&$9$&$10$&$11$&$12$&$13$&$14$&$15$&$16$\\
\hline
&$7$&$0$&$16$&$19$&$32$&$35$&$48$&$51$&$64$&$67$&$80$&$83$&$96$&$99$&$112$&$115$&$128$\\
&$8$&$0$&$15$&$18$&$32$&$34$&$47$&$50$&$64$&$66$&$79$&$82$&$96$&$98$&$111$&$114$&$128$\\
$e$&$9$&$0$&$14$&$17$&$31$&$33$&$46$&$49$&$64$&$66$&$78$&$81$&$95$&$97$&$110$&$113$&$128$\\
&$10$&$0$&$13$&$16$&$30$&$32$&$45$&$48$&$63$&$65$&$77$&$80$&$94$&$96$&$109$&$112$&$128$\\
&$11$&$0$&$12$&$16$&$29$&$31$&$44$&$47$&$62$&$64$&$76$&$79$&$93$&$95$&$108$&$111$&$127$\\
&$12$&$0$&$12$&$16$&$28$&$30$&$43$&$46$&$61$&$63$&$75$&$78$&$92$&$94$&$107$&$110$&$126$\\
&$13$&$0$&$12$&$16$&$27$&$29$&$42$&$45$&$60$&$62$&$74$&$77$&$91$&$93$&$106$&$109$&$125$\\
&$14$&$0$&$12$&$16$&$26$&$28$&$41$&$44$&$59$&$61$&$73$&$76$&$90$&$92$&$105$&$108$&$124$\\
\hline
&&&&&&$k$&&&&&&&&&&&\\
&&$17$&$18$&$19$&$20$&$21$&$22$&$23$&$24$&$25$&$26$&$27$&$28$&$29$&$30$&$31$&$32$\\
\hline
&$8$&$130$&$143$&$146$&$160$&$162$&$175$&$178$&$192$&$194$&$207$&$210$&$224$&$226$&$239$&$242$&$256$\\
&$9$&$130$&$142$&$145$&$159$&$161$&$174$&$177$&$192$&$194$&$206$&$209$&$223$&$225$&$238$&$241$&$256$\\
$e$&$10$&$130$&$141$&$144$&$158$&$160$&$173$&$176$&$191$&$193$&$205$&$208$&$222$&$224$&$237$&$240$&$256$\\
&$11$&$129$&$140$&$143$&$157$&$159$&$172$&$175$&$190$&$192$&$204$&$207$&$221$&$223$&$236$&$239$&$256$\\
&$12$&$128$&$139$&$142$&$156$&$158$&$171$&$174$&$189$&$191$&$203$&$206$&$220$&$222$&$235$&$238$&$255$\\
&$13$&$127$&$138$&$141$&$155$&$157$&$170$&$173$&$188$&$190$&$202$&$205$&$219$&$221$&$234$&$237$&$254$\\
&$14$&$126$&$137$&$140$&$154$&$156$&$169$&$172$&$187$&$189$&$201$&$204$&$218$&$220$&$233$&$236$&$253$
\end{tabular}
}
\caption{Values of $g(e,k)$ when $e\le14$ and $k\le32$.}
\label{gtbl}
\end{center}
\end{table}

In Section \fullref{normsec}, we apply the basic induction argument, \fullref{compatprop}, and the results for $\gd(2^e\xi)$ in \fullref{gthm}
to prove the following result by induction on $t$. This clearly implies \fullref{gdthm} and hence \fullref{mainthm}.

\begin{prop}\label{1better} For $e\ge 7$ and $t\ge0$,
$\gd((2^e+2^{e+1}+\cdots+2^{e+t})\xi_{2^e-1})\le 2^e-e-6$.
\end{prop}

\section{Proof of general lifting results}\label{compatsec}
In this section, we prove \fullref{compatprop} and \fullref{pairedprop}. For the first one, we find it more convenient to
work with sections rather than geometric dimension.
\begin{thm} \label{sectthm}Let $X_0\subset\cdots\subset X_k$ and
$Y_0\subset\cdots\subset Y_k$
be filtered spaces, and let $\theta$ (resp. $\eta$) be a vector bundle over
$X_k$ (resp. $Y_k$). Suppose given $m_0$ (resp.~$n_0$) sections of
$\theta$ on $X_k$ (resp. $\eta$ on $Y_k$), of which the first $m_i$
(resp. $n_i$) are linearly independent (l.i.) on $X_{i}$ (resp. $Y_{i}$)
for
$0\le i\le k$. Let $$p_j=\min(m_{i}+n_{j-i}: 0\le i\le j).$$ Let
$$W_j=\bigcup_{i=0}^j X_i\times Y_{j-i}.$$ Then there are $p_0$ sections of
$\theta\times\eta$ on $W_k$ of which the first $p_j$ are linearly
independent on $W_{j}$ for $0\le j\le k$. Moreover, if $\ell+i\ge j$ and
$m_\ell+n_i\ge p_j$, then the first $p_j$ sections are l.i.~on
$X_\ell\times Y_i$.
\end{thm}

Note that we have $m_0\ge \cdots\ge m_k$, $n_0\ge \cdots\ge n_k$, and
$p_0\ge\cdots\ge p_k$.

The following result will be used in the final step of the proof of
\fullref{sectthm}.
\begin{lem} Suppose  $\theta$ is an $n$--dimensional
trivial vector bundle over a space $X$ with l.i.~sections $t_1,\ldots,t_n$.
Suppose $s_1,\ldots,s_r$ are l.i.~sections of $\theta$, each of
which is a linear combination with constant coefficients of the $t_i$. Then
there is a set $s_1,\ldots,s_r,s'_{r+1},\ldots,s'_n$ of linearly
independent sections of $\theta$, with all these sections being linear
combinations with constant coefficients of the $t_i$.\end{lem}

\begin{proof} Because of the constant-coefficient assumption, this is just
a consequence of
the result for vector spaces, that a basis for a subspace can be extended
to a basis for the
whole space.
\end{proof}

Note that the assumption about constant coefficients was required. For
example, the section $s(x)=(x,x)$
of $S\sp2\times \R\sp3$ cannot be extended to a set of three
l.i.~sections.

\begin{proof}[Proof of \fullref{sectthm}] Let $r_1,\ldots,r_{m_0}$ be the
given sections of $\theta$ on $X_k$,
and $s_1,\ldots,s_{n_0}$ the given sections of $\eta$ on $Y_k$. These are
considered as sections of $\theta\times\eta$ by using 0 on the other
component. Clearly $\{r_1,\ldots,r_{m_0},s_1,\ldots,s_{n_0}\}$ is a set of
$p_0$ sections on $W_k$ which is linearly independent on $W_0$. The
proof will proceed by finding $p_1$ linear combinations, always with
constant coefficients, of these sections which are l.i.~on $W_1$, then
$p_2$ linear combinations of these new sections which are l.i.~on $W_2$,
etc, until going into the last stage we have $p_{k-1}$ sections which
are l.i.~on $W_{k-1}$, and we find $p_k$ linear combinations of them which
are l.i.~on $W_k$. Now we apply the lemma repeatedly, starting with
the last $p_k$ sections. At the first step, we extend this set to a set of
$p_{k-1}$ sections l.i.~on $W_{k-1}$, and continue until going into
the last stage, where we have $p_1$ sections which are combinations of the
original $p_0$ sections and satisfy the conclusion of the theorem for
$1\le i\le k$. We apply the lemma one last time to extend the set of $p_1$
sections to the desired set of $p_0$ sections.

Here is an explicit algorithm for the sections described in the first half
of the preceding paragraph.
We may assume without loss of generality that $m_0\ge n_0$.

For $j$ from 0 to $k$,
\begin{itemize}
\item For $i$ from 1 to $p_j-n_0$ (resp. $p_j-m_0$), let $r_i^{(j)}=r_i$
(resp. $s_i^{(j)}=s_i$). (Note that if $n_0\ge p_j$,
then nothing happens at this step.)
\item For $i$ from $\max(1,p_j-n_0+1)$ to $\min(m_0,p_j)$, let both
$r_i^{(j)}$ and $s_{p_j+1-i}^{(j)}$
equal $r_i^{(j-1)}+s_{p_j+1-i}^{(j-1)}$.
\item Then the sections $r_i^{(j)}$ and $s^{(j)}_i$ constructed in the two
previous steps give the sections
which are l.i.~on $W_j$. (Each section constructed in the second step can
be counted as an $r$ or an $s$,
but is only counted once.)
\end{itemize}

We must show that these have the required linear independence. Before doing
so, we illustrate with
an example, computed by {\tt Maple}. Let $k=4$,
$[m_0,\ldots,m_4]=[11,6,4,1,0]$ and $[n_0,\ldots,n_4]=[10,8,3,2,0]$. Then
$[p_0,\ldots,p_4]=[21,16,14,9,7]$.
The 16 sections l.i.~on $W_1$ are
$$r_1,\ldots,r_6,\ r_7+s_{10},\ r_8+s_{9},\ r_9+s_8,\ r_{10}+s_7,\
r_{11}+s_6,\ s_5,\ldots,s_1.$$
The 14 sections l.i.~on $W_2$ are
\begin{eqnarray*}&&r_1,\ r_2,\ r_3,\ r_4,\ r_5+r_7+s_{10},\ r_6+r_8+s_{9},\
r_7+r_9+s_{10}+s_8,\ r_8+r_{10}+s_9+s_7,\\
&&r_9+r_{11}+s_8+s_6,\ r_{10}+s_7+s_5,\ r_{11}+s_6+s_4,\ s_3,\ s_2,\
s_1.\end{eqnarray*}
The 9 sections l.i.~on $W_3$ are
\begin{eqnarray*}&&r_1+r_6+r_8+s_{9},\  r_2+r_7+r_9+s_{10}+s_8,\
r_3+r_8+r_{10}+s_{9}+s_7,\\
&&r_4+r_9+r_{11}+s_8+s_6,\ r_5+r_7+r_{10}+s_{10}+s_7+s_5,\\
&& r_6+r_8+r_{11}+s_9+s_6+s_4,
r_7+r_9+s_{10}+s_8+s_3,\ r_8+r_{10}+s_9+s_7+s_2,\\\
&&r_9+r_{11}+s_8+s_6+s_1.\end{eqnarray*}
The 7 sections l.i.~on $W_4$ are
\begin{eqnarray*}&&r_1+r_3+r_6+2r_8+r_{10}+2s_9+s_7,\\
&&r_2+r_4+r_7+2r_9+r_{11}+s_{10}+2s_8+s_6,\\
&&r_3+r_5+r_7+r_8+2r_{10}+s_{10}+s_9+2s_7+s_5,\\
&&r_4+r_6+r_8+r_9+2r_{11}+s_{9}+s_8+2s_6+s_4,\\
&&r_5+2r_7+r_9+r_{10}+2s_{10}+s_8+s_7+s_5+s_3,\\
&&r_6+2r_8+r_{10}+r_{11}+2s_9+s_7+s_6+s_4+s_2,\\
&&r_7+2r_9+r_{11}+s_{10}+2s_8+s_6+s_3+s_1.\end{eqnarray*}

Now we continue with the proof. The property described in the first
paragraph of the proof, that
the sections claimed to be l.i.~on $W_j$ are linear combinations with constant coefficients of those
on $W_{j-1}$, is clear from
their inductive definition.

Next we easily show that if $i>p_j-n_0$, then
$$r_i^{(j)}=s_{p_j+1-i}^{(j)}=r_i+\sum_{\ell>i}c_\ell
r_\ell+s_{p_j+1-i}+\sum_{\ell>p_j+1-i}d_\ell s_\ell$$
with $c_\ell$ and $d_\ell$ integers. The point here is that the additional
terms have subscript
greater than $i$ or $p_j+1-i$. The proof is immediate from the inductive
formula
$$r_i^{(j)}=r_i^{(j-1)}+s_{p_j+1-i}^{(j-1)}$$
and the fact that $p_j\le p_{j-1}$. Indeed, from $r_i^{(j-1)}$ we obtain
terms $r_{\ge i}$ and
$s_{\ge p_{j-1}+1-i}$, and from $s_{p_j+1-i}^{(j-1)}$ we obtain terms
$s_{\ge p_j+1-i}$ and
$r_{\ge p_{j-1}-p_j+i}$.

Finally we show that the asserted sections are l.i.~on $W_j$. Let $\bx\in
X_\ell\times Y_{j-\ell}$.
Note that $\{r_1(\bx),\ldots,r_{m_\ell}(\bx)\}$ is l.i., as is
$\{s_1(\bx),\ldots,s_{n_{j-\ell}}(\bx)\}$,
and that $p_j\le m_\ell+n_{j-\ell}$. If we form a matrix with columns
labeled
$$r_1,\ldots,r_{m_0},s_{n_0},\ldots,s_1,$$
and rows which express the sections, ordered as
\begin{equation}\label{secor}r_1^{(j)},
\ldots,r_{\min(m_0,p_j)}^{(j)},s^{(j)}_{p_j-m_0},\ldots,s_1^{(j)},\end{equation}
in terms
of the column labels, then, by the previous paragraph, the number of
columns is $\ge$ (usually strictly greater than) the number
of rows, the entry in position $(i,i)$ is 1 for $i\le\min(m_0,p_j)$, and
all entries to the left of these
1s are zero. If $i>\min(m_0,p_j)$, then all entries in the $r$-portion of row
$i$ are zero.
Moreover an analogous statement is true if the order of the rows and of the
columns are both reversed.
Thus there are 1s on the diagonal running
up from the lower right corner of the original matrix (for $\min(n_0,p_j)$
positions) and zeros to their right.

If a linear combination of our sections applied to $\bx$ is 0, then the
triangular form of the matrix
implies that the first $m_\ell$ coefficients are 0, while the triangular
form looking up from the lower
right corner implies that the last $n_{j-\ell}$ coefficients are 0. Since
$p_j\le m_\ell+n_{j-\ell}$, this implies
that all coefficients are 0, hence the desired independence.

The same argument works for the last statement of the proposition. For $k$
satisfying $j\le k\le\ell+i$, replace
$W_k$ by $W_k\cup (X_{\ell}\times Y_i)$. Then everything goes through as
above.
\end{proof}

\begin{proof}[Proof of \fullref{compatprop}] Let $D=\dim(\theta)$ and
$D'=\dim(\eta)$. Then
$d_i$, $d'_i$, $e_i$, and $(X\times Y)_{i}$ of \fullref{compatprop}
correspond to  $D-m_i$, $D'-n_i$, $D+D'-p_i$, and $W_i$ of
\fullref{sectthm}, respectively. The compatible gd bounds may be
interpreted
as vector bundles $\theta_i$ over $X_i$ of dimension $d_i$ and isomorphisms
$\theta|X_i\approx \theta_i\oplus(D-d_i)$
and $\theta_{i}|X_{i-1}\approx\theta_{i-1}\oplus(d_{i}-d_{i-1})$. The
trivial subbundles yield, for all $i$, $D-d_i$
l.i.~sections of $\theta$ on $X_i$ such that the restrictions of the
sections on $X_i$ to $X_{i-1}$ are a subset of the sections on $X_{i-1}$.
Each of the sections on $X_0$ has a largest $X_i$ for which it is one of
the given l.i.~sections. By Atiyah
\cite[Section 1.4.1]{At}, this section on $X_i$ can be
extended over $X_k$ (although probably not as part of a linearly
independent set). Analogous statements are true for sections of $\eta|Y_i$.

By \fullref{sectthm}, there are $D+D'-e_0$ l.i.~sections of
$\theta\times\eta$ on $W_0$ of which the first $D+D'-e_i$ are l.i.~on
$W_i$. Taking orthogonal complements of the spans of the sections yields
the desired
compatible bundles on $W_i$ of dimension $e_i$, yielding the first part of
\fullref{compatprop}.

For the second part, first note that in the algorithm in the proof of
\fullref{sectthm}, if the $r$'s and $s$'s are equal, then the set of
sections constructed on each $W_i$ is invariant under the interchange map
$T$. Thus the same will be true of the orthogonal complement of their
span.\end{proof}

\begin{proof}[Proof of \fullref{pairedprop}] \qua
\begin{enumerate}
\item Let $F_1=S^{2m-1}$ denote the fiber of $BO_{2m-1}[\rho]\to
BO_{2m}[\rho]$. There is a relative
Serre spectral sequence for
\begin{equation}\label{SSS1}(CF_1,F_1)\to
(BO_{2m}[\rho],BO_{2m-1}[\rho])\to BO_{2m}[\rho].\end{equation}
The fibration $V_{2m}\to BO_{2m}[\rho]\to BO[\rho]$ shows that the bottom
class of $BO_{2m}[\rho]$ is in dimension $\ge\min(\rho,2m)$. The
spectral sequence of \eqref{SSS1} shows that $H_*(S^{2m})\to
H_*(BO_{2m}[\rho]/BO_{2m-1}[\rho])$ has cokernel beginning in dimension
$\ge s+1$,
and so the map is an $s$--equivalence. Thus the inclusion of the
$s$--skeleton of $BO_{2m}[\rho]/BO_{2m-1}[\rho]$ factors through $S^{2m}$ to
yield the map $c_1'$, which is an equivalence.

The second map is obtained similarly. A map $$\Sigma
P_{2m-1}^{2m}\mapright{\ell}BO_{2m+1}[\rho]/BO_{2m-1}[\rho]$$ is obtained as
the inclusion of
a skeleton of $CF_2/F_2$, where $F_2=V_{2m+1,2}$ is the fiber of
$BO_{2m-1}[\rho]\to BO_{2m+1}[\rho]$. The relative Serre spectral sequence
of
\begin{equation}(CF_2,F_2)\to (BO_{2m+1}[\rho],BO_{2m-1}[\rho]) \to
BO_{2m+1}[\rho]\label{SSS}\end{equation}
implies that $\coker(\ell_*)$ begins in dimension $\ge s+1$, determined by
$$H_{2m}(CF_2,F_2)\otimes H_{\min(\rho,2m+1)}(BO_{2m+1}[\rho])$$  and
the first ``product" class in $H_{4m}(\Sigma V_{2m+1,2})$. The obtaining of
$c_2'$ now follows exactly as for $c_1'$.

\item Let $Q :=BO_{2m+1}[\rho]/BO_{2m-1}[\rho]$ and
$E:=\text{fiber}(BO_{2m+1}[\rho]\to Q)$.
The commutative diagram of fibrations
$$\begin{CD} V_{2m+1,2}@>>> BO_{2m-1}[\rho] @>>> BO_{2m+1}[\rho]\\
@VVV @VVV @VVV\\
\Omega Q@>>> E @>>> BO_{2m+1}[\rho]\end{CD}$$
implies the quotient $E/BO_{2m-1}[\rho]$ has the same connectivity as
$\Omega Q/V_{2m+1,2}$, which is 1 less than that
determined from \eqref{SSS}; that is, $E/BO_{2m-1}[\rho]$ is
$(s-1)$--connected. Thus, since $\dim(X)<s$, the vertical
maps in
$$\begin{CD}BO_{2m-1}[\rho]^{(s)}@>>>BO_{2m+1}[\rho]^{(s)}@>>>\Sigma
P_{2m-1}^{2m}\\
@VVV @VVV @VVV\\
E @>>> BO_{2m+1}[\rho] @>>> Q\end{CD}$$
are equivalences in the range relevant for maps from $X$, $A$, and $X/A$.
Since the bottom row is a fibration,
we may consider the top row to be one, too, as far as $X$ is concerned.

Since $g$ is divisible by 2, and $2\pi_{2m}(\Sigma P_{2m-1}^{2m})=0$, we
deduce that the composite
$$X/A\mapright{g} S^{2m} \mapright{i} \Sigma P_{2m-1}^{2m}$$
represents the 0 element of $[X/A,\Sigma P_{2m-1}^{2m}]$; ie the map is
null-homotopic rel $*$. There is a commutative diagram as in
\eqref{coffib} with the left sequence a cofiber sequence and the right
sequence a fiber sequence in the range of $\dim(X)$.
\begin{equation}\label{coffib}\begin{CD}
A @>f_1>> BO_{2m-1}[\rho]^{(s)}\\
@Vj_1VV @Vj_2VV\\
X @>p_2\circ f>> BO_{2m+1}[\rho]^{(s)}\\
@VqVV @VVV \\
X/A@>i\circ g>> \Sigma P_{2m-1}^{2m}\end{CD}\end{equation}
We have just seen that there is a basepoint-preserving homotopy
$$H \co X/A\times I\to \Sigma P_{2m-1}^{2m}$$
from $i\circ g$ to a constant map.
There is a commutative diagram
$$\begin{CD} X\times 0\cup A\times I@.\longrightarrow @.BO_{2m+1}[\rho]\\
@VVV @. @VVV\\
X\times I @>q\times I>> X/A\times I @>H>> \Sigma P_{2m-1}^{2m}\end{CD}$$ where the top map is $p_2\circ f$ on $X\times 0$ and $j_2\circ f_1$ on
each $A\times\{t\}$. By the Relative Homotopy Lifting Property of a fibration, there exists a map $\widetilde H \co X\times I\to
BO_{2m+1}[\rho]$ making both triangles commute. When $t=1$, it maps into $BO_{2m-1}[\rho]$, since it projects to the constant map at the
basepoint of $\Sigma P_{2m-1}^{2m}$.

\item
We use the fact that $2\cdot 1_{\Sigma P_{2m-1}^{2m}}$
factors as
$$\Sigma P_{2m-1}^{2m}\mapright{\text{col}}S^{2m+1}\mapright{\eta}
S^{2m}\hookrightarrow \Sigma P_{2m-1}^{2m}$$
to deduce that the composite
$$\Sigma X/A\mapright{\Sigma g'}\Sigma^2 P_{2m-1}^{2m}\mapright{\text{col}}
S^{2m+2}$$
is null-homotopic since $[\Sigma g']$ is divisible by 2. Note that we
needed to suspend once since if $\dim(X)=4m-1$, then $[X/A,\Sigma
P_{2m-1}^{2m}]$ might not have a group structure. Since
$$[X/A,S^{2m+1}]\mapright{\Sigma}[\Sigma (X/A),S^{2m+2}]$$
is bijective, we deduce that $X/A\mapright{{\text{col}\circ g'}} S^{2m+1}$
is null-homotopic.

An argument similar to the one in the beginning of the proof of (2) shows that $BO_{2m}[\rho]\to BO_{2m+1}[\rho]\to S^{2m+1}$ is a fibration
through dimension $\min(\rho+2m-1,4m)\ge s$. Since $\dim(X)\le s$, the lifting follows as in the proof of (2). \proved
\end{enumerate}
\end{proof}

\section[Inductive determination of a bound for gd(2^e,k)]{Inductive
determination of a bound for $\gd(2^e,k)$}\label{indnsec}

In this section, we prove that part (6) of \fullref{gthm} follows from its first five parts, together with its validity for $e=7$.
We begin by proving the validity when $e=7$.
The following result is stronger than the required liftings for $e=7$; i.e., we have $m(k)\le g(7,k)$ and the inequality is strict
if $k$ is even with $4\le k\le 14$. The reason for beginning our induction with liftings weaker than the best results that we are able to prove
is to fit them into a simple formula that works for all values of $e$. Here and throughout we use the standard notation that $\nu(-)$ denotes
the exponent of 2 in an integer.

\begin{thm}\label{128prop} Let
$$m(k)=\begin{cases}0&k=0,1\\
16&k=2\\
8k-5&k\text{ odd, }\ 3\le k\le15\\
8k+\nu(k)-4&k\text{ even, }\ 4\le k\le16.\end{cases}$$
There are compatible liftings of $128\xi_{8k+7}$ to $BO_{m(k)}$ for $k\ge0$.\end{thm}
\begin{proof}
Let $H_k$ denote the Hopf bundle over quaternionic projective space $HP^k$.
Let $m'(k)=13$ if $k=2$, and otherwise $m'(k)=m(k)$.
We will use \cite[Theorem 1.1b]{HP} to prove
\begin{equation}\label{lifts}\text {there are compatible liftings of $32H_{2k+1}$ to $BO_{m'(k)}$ for $2\le k\le 16$.}\end{equation}

Three things are required to prove this. First we need that, for $k\le 15$ and  all $i\le 2k+1$
satisfying also $4i-1\ge m'(k)$,
$$\nu\tbinom{32}i\ge\nu(|\pi_{4i-1}(P_{m'(k)}\w bo)|).$$
This is easily verified using $\nu\tbinom{32}i=5-\nu(i)$ and, for $1\le\eps\le3$,
\begin{equation}\nu(|\pi_{4i-1}(P_{4a+\eps}\w bo)|)=\begin{cases}4-\eps&i=a+1\\
4&i=a+2\\ 8-\eps&i=a+3.\end{cases}\label{hts}\end{equation} For example, if $k$ is odd, we have $m'(k)=8k-5$. Then $a=2k-2$ and $\eps=3$ in
(\ref{hts}), and for $i=\langle 2k-1,2k,2k+1\rangle$, we have $\nu\tbinom{32}i=\langle 5,4,5\rangle$ and $\nu(|\pi_{4i-1}(P_{8k-5}\w
bo)|)=\langle 1,4,5\rangle$.

Secondly, we need that $\pi_{4i-1}(P_{m'(k)})\to\pi_{4i-1}(P_{m'(k)}\w bo)$ is injective for $i\le 2k+1$. This is obtained from Tables 8.4, 8.8,
8.14, 8.15, and 8.16 of \cite{Mem}. These show that for $m'(k)\equiv\langle3,7,13,14,15\rangle$ mod 16 and $4i-1\le m'(k)+\langle8,4,6,5,4\rangle$, the
asserted injectivity is true. Now the liftings follow from \cite[Theorem 1.1b]{HP}. If $k=16$, the lifting follows for dimensional reasons.

The third thing we need is compatibility. We must show that
$$\begin{CD}HP^{2k-1}@>>> BO_{m'(k-1)}\\
@VVV @VVV\\
HP^{2k+1}@>>> BO_{m'(k)}\end{CD}$$ commutes for $k\ge3$. The two composites agree stably, and so their obstructions to being homotopic lie in
$H^*(HP^{2k-1};\pi_*(V_{m'(k)}))$. If $k$ is even, then $8k-4<m'(k)$ so the groups are 0. If $k$ is odd, the result follows since
$\pi_{8k-4}(V_{8k-5})=0$.

We precede the compatible liftings of (\ref{lifts}) by the canonical maps $RP^{8k+7}\to HP^{2k+1}$, obtaining compatible liftings of
$128\xi_{8k+7}$ to $BO_{m(k)}$ for $k\ge2$. The bundle $128\xi_{15}$ is trivial. To insure compatibility of the liftings on $RP^{15}$ and $RP^{23}$,
we note that the obstructions to compatibility lie in $H^*(RP^{15};\pi_*(V_{16}))=0$.
This is why we use $m(k)=16$, rather than 13.
\end{proof}

Now we prove the induction step. Let
$$\rho(4a+b)=8a+2^b\text{ if }0\le b\le3.$$
It satisfies that $2^k\xi_n$ is nontrivial if and only if $n\ge\rho(k)$.
Let $\rho=\rho(e-1)$. Assume that we have obtained compatible liftings of $2^{e-1}\xi_{8k+7}$ to
$BO_{g(e-1,k)}[\rho]$ for all $k$. For $0\le k\le 2^{e-3}$, define $$g_1(e,k):=\max\{g(e-1,i)+g(e-1,k-i): \max(0,k-2^{e-4})\le i\le[k/2]\}.$$
Note that by \fullref{gthm}.(3),
\begin{equation}\label{gg1}g(e,k)\ge g_1(e,k)-1.\end{equation}

Recall $A_k=P^{8k+7}$, and let
$$(A\times A)_k=\bigcup_{i=0}^k A_i\times A_{k-i}.$$
Then by \fullref{compatprop} there are compatible symmetric liftings
$\ell_k$ of $2^{e-1}\xi\times 2^{e-1}\xi$ on $(A\times A)_{k}$ to
$BO_{g_1(e,k)}[\rho]$ for all $k$. We precede by compatible maps
$d_k:A_k\to (A\times A)_{k}$, cellular maps homotopic to the diagonal. The
composites
$A_k\stackrel{\ell_k\circ d_k}{\longrightarrow} BO_{g_1(e,k)}[\rho]$
are
compatible liftings of $2^e\xi_{8k+7}$ for all $k$.

By decreasing induction on $k$ starting with $k=2^{e-3}$, we will construct
compatible factorizations through
$BO_{g(e,k)}[\rho]$ of the maps $\ell_k\circ d_k$.
Assume inductively that, for all $j>k$, compatible factorizations, up to
homotopy rel $A_k$, of $\ell_j\circ d_j$
through $BO_{g(e,j)}[\rho]$ have been attained. If $g(e,k)\ge g_1(e,k)$,
then no factorization of $\ell_k\circ d_k$ is required, and so our
induction on $k$ is extended. So we may assume $g(e,k)=g_1(e,k)-1$.

Let $h=[k/2]$.  By \fullref{gthm}.(4),
\begin{equation}\label{2diff}g_1(e,k-1)\le g(e,k)-1.\end{equation}
By \eqref{2diff}, \fullref{gthm}.(4), and the last part of \fullref{sectthm}
(which is required for compatibility of the
lifts of $(A\times A)_{k-1}$ and $A_h\times A_h$ to $BO_{g(e,k)-1}$), we
have the commutative diagram below, similar to \eqref{coffib}.
$$\begin{CD}A_{k-1}@>d'>>(A\times A)_{k-1}\cup A_h\times
A_h@>>>BO_{g(e,k)-1}[\rho]^{(8k+7)}\\
@VVV @VVV @VVV\\
A_k @. (A\times A)_k @. BO_{g(e,k)}[\rho]^{(8k+7)}\\
@| @| @VVV\\
A_k@>d_k>>(A\times A)_{k}@>\ell_k>> BO_{g(e,k)+1}[\rho]^{(8k+7)}\\
@VVV @VVV @VcVV\\
A_k/A_{k-1}@>\dbar>>(A\times A)_{k}/((A\times A)_{k-1}\cup A_h\times A_h)@>\lbar>> C,\end{CD}$$ where $C=S^{g(e,k)+1}$ if $g(e,k)$ is odd, and
$C=\Sigma P^{g(e,k)}_{g(e,k)-1}$ if $g(e,k)$ is even. The maps labeled $d$ are cellular maps homotopic to the diagonal. The map $c$ is obtained
similarly to the first paragraph of the proof of \fullref{pairedprop}. Since $\dim(A_k)=8k+7$, the application of \fullref{pairedprop}
requires that
$$8k+7\le\min(\rho+g(e,k)-1,2g(e,k)-1).$$
The second follows from \fullref{gthm}.(2), while the first follows from $\rho\ge 2e-2$ and $g(e,k)\ge 8k-e+2$ since $e\ge8$.

The quotient $(A\times A)_{k}/(A_h\times A_h)$ equals $B\vee T(B)$,
where $T$ reverses the order of the factors, and
$B$ is the union of all cells $e^i\times e^j$ with $i<j$.
By the symmetry property of $\ell_k$, $\lbar|T(B)=(\lbar|B)\circ T$. Since
$T\circ\dbar \simeq \dbar$, we conclude that
$\lbar\circ\dbar$ is divisible by 2. Indeed, with $r_B$ denoting the
retraction onto $B$,
$$[\lbar\circ\dbar]=[(\lbar|B)\circ r_B\circ\dbar]+[(\lbar|T(B))\circ
r_{T(B)}\circ\dbar]$$
and we have $$[(\lbar|T(B))\circ r_{T(B)}\circ\dbar]=[(\lbar|T(B))\circ
T\circ r_B\circ \dbar]=[(\lbar|B)\circ r_B\circ\dbar].$$

Thus, by \fullref{pairedprop}, $\ell_k\circ d_k$ is homotopic rel $A_{k-1}$
to a map which lifts to $BO_{g(e,k)}[\rho]$. Note that the
lifting into $BO_{g(e,k)-1}[\rho]$ was not needed if $g(e,k)$ is odd. We have extended our
inductive lifting hypothesis, and so  have proved that there are
compatible liftings of $A_k$ to $BO_{g(e,k)}[\rho]$ for all $k$. This
extends the induction on $e$ and proves \fullref{gthm}.(6), assuming the
first five parts of \fullref{gthm}.

\section{The function $g(e,k)$}\label{numbersec}

In this section, we define the function $g(e,k)$ which has been used in the
previous
sections, and prove the first five parts of \fullref{gthm}, its numerical
properties
which were already used to prove \fullref{gthm}.(6), its important geometrical
property.

We find it convenient to deal with the complementary function $G$ defined by
\begin{equation}\label{complement}G(e,k)=8k-g(e,k).\end{equation}
It has relatively small values, in which patterns are more readily apparent.
This function $G$ will be defined using several auxiliary functions.

We define a function $S$ for $k\ge2$ by
\begin{equation}\label{Sdef}S(k)=8k-13[\tfrac{k+1}2]+2\alpha(k)+2\min(3,\nu(k-1))+\begin{cases} -1&k\equiv0\ (2)\\2&k\equiv1\ (8)\text{ and } \alpha(k)\ne2\\
4&\text{otherwise.}\end{cases}\end{equation}
Then $S(k)=8k-s(k)$, where $s(k)$ is the stable value of $g(e,k)$ when $e$ is sufficiently large.
The first values of $S$ are given by

\begin{center}
\begin{tabular}{c|cccccccccccccccccc|}
\hline
$k$&2&3&4&5&6&7&8&9&10&11&12&13&14&15&16&17&18&19\\
\hline
$S(k)$&4&8&7&13&12&16&13&21&18&22&21&27&26&30&25&33&30&34\\
\hline
\end{tabular}
\end{center}
It will occasionally be useful to set $S(1)=8$, consistent with $g(e,1)=0$.

Values of $k\equiv1$ mod 8 receive special treatment. They are excluded in the domain of some of
our functions. For example, for $k\not\equiv 1$ mod 8 with $k\ge2$, we define $V(k)$ by
$$V(k)=\a(k)-\begin{cases}2&k\equiv3\ (4)\\ 1&k\not\equiv3\ (4).\end{cases}$$
The reasons for defining some of these functions will be presented shortly.

We also define functions $\nu'$ and $R$ as follows.
$$\nu'(k)=\begin{cases}\nu(k)&k\text{ even}\\
-4&k\text{ odd,}\end{cases}$$
and, for $k\not\equiv1$ mod 8, \begin{equation}\label{Rdef}R(k)=S(k)+\nu'(k)-V(k).\end{equation}
The first few values of $R$ are given as follows.

\begin{center}
\begin{tabular}{c|cccccccccccccccc|}
\hline
$k$&2&3&4&5&6&7&8&10&11&12&13&14&15&16&18&19\\
\hline
$R(k)$&5&4&9&8&12&11&16&18&17&22&21&25&24&29&30&29\\
\hline
\end{tabular}
\end{center}

It will also be useful to introduce the notation $\la n\ra=\max(0,n)$.
We will frequently use the simple fact that for any number $X$,
\begin{equation}\label{X}X+\la -X\ra\ge0.\end{equation}

Now we can define our function $G$.
An integer $k$ is {\bf decomposable} if it can be written as $k=k_0+\cdots+k_r$ with $r\ge1$ and $\nu(k_i)>R(k_{i-1})$ for $1\le i\le r$.
Because each $k_i$  must be preceded in the binary expansion of $k$ by a long string of 0's,
it is clear that a decomposable integer has a unique maximal decomposition. The sum in (\ref{geq1}) is taken over all $k_i$, $i\ge1$,
in a maximal decomposition of $k$.
The first values of $k$ admitting a decomposition are 35, 66, and 67, with $k_0=3$, 2, 3,
respectively. A simple decomposition is a maximal decomposition with $r=1$.
 The first value of $k$ admitting a multiple decomposition is $2^{55}+35$ with $k_0=3$ and $k_1=32$.

\begin{defin}\label{geq}
If $2\le k\le 2^{e-3}$ and $(e,k)\ne(7,9)$, we define
$$G(e,k)=\begin{cases}\min(S(k),G'(e,k))&k\not\equiv1\ (8)\\
\min(S(k),6+G'(e,k-1))&k\equiv1\ (8),\end{cases}$$
where, for $k\not\equiv1$ mod 8,
\begin{equation}G'(e,k)=\langle e-6-\nu'(k)\rangle-\sum_{k_i}\langle\min(\nu(k_i),e-6)-R(k_{i-1})\rangle.\label{geq1}\end{equation}
The exceptional value is $G(7,9)=5$, not 6 as the formula would give.
\end{defin}

The terms in the sum in \eqref{geq1} will sometimes be called {\it deviations}.
We do not define $G(e,0)$, $G(e,1)$, or $G(e,k)$ for $k>2^{e-3}$; instead we just define the complementary function $g$ by $g(e,0)=g(e,1)=0$
and $g(e,k)=2^{e-3}$ for $k>2^{e-3}$, and observe that the crucial properties (3) and (4) in Theorem \ref{gthm} are easily seen to be
satisfied whenever these extreme values are involved.

Next we provide some general discussion of what led to the rather complicated formula for $g(e,k)$. First we describe what led to the
basic formula $g(e,k)\approx 8k-\la e-6-\nu'(k)\rangle$, modified when $k\equiv1$ mod 8.
We began with the initial values $m(k)$ of \fullref{128prop} for $g(7,k)$ and used a computer
program implementing properties (3) and (4) of \fullref{gthm} to obtain bounds for $g(e,k)$ for larger $e$.
Except perhaps for the first few entries in a $k$-column, the values $8k-(e-6-\nu'(k))$ when $k\not\equiv1$ mod 8, and $g(e,8\ell+1)=g(e,8\ell)+2$,
were apparent until issues of stabilization, which we will discuss shortly, became involved. However, there was no apparent regular pattern
for the first few entries in each $k$-column. The formula $8k-\la e-6-\nu'(k)\ra$ was achieved after additional computer experimentation
as the simplest general formula satisfying $g(7,k)\ge m(k)$ and consistency with \fullref{gthm}.

Next we explain where $S(k)$ came from. It is related to the condition $g(e,k)\ge 4k+4$, which says that our lifting methods only work
in the stable range. In an earlier version of this paper, we used the triviality of $2^{\phi(n)}\xi_n$ to give 0 as the value of
$g(e,k)$ when $e>4k+3$, but we were unable to prove that this could be done compatibly with our other liftings; i.e. that the liftings
which we obtain inductively can be done so that their restrictions to appropriate skeleta are trivial. By forcing $g(7,2)=16$, we could,
as noted in the proof of \fullref{128prop}, guarantee that our liftings restrict to a trivial bundle on $P^{15}$, the case $k=1$.
For reasons of stability, we forced $g(e,2)\ge 12$ and $g(e,3)\ge 16$. Forcing $g(e,4)\ge20$ is not strong enough, since, with $g(15,2)=12$
and $g(15,4)=25$, we could not obtain $g(16,4)=24$ consistently with property (4) of \fullref{gthm}. Thus $g(e,4)=25$ for $e\ge15$; i.e. $s(4)=25$.
This translates to our value $S(4)=8\cdot 4-s(4)=7$.

To be consistent with \fullref{gthm}, our function $S$ must satisfy the inequalities of the following proposition, the proof of which
is straightforward, although somewhat tedious, and is omitted.
\begin{prop}\label{Sprop} The function $S$ defined in \eqref{Sdef} satisfies
$$S(i+j)\le S(i)+S(j)+1$$
and $$S(2i)\le 2S(i)-1,$$
with equality in the first if $i=2^t$ and $2\le j\le 2^t-1$ or $j=2^t+1$, and equality holds in the second if $i=2^t$.
Thus $S$ may be defined by $S(2)=4$, $S(3)=8$, and
$$S(k)=\min(S(i)+S(k-i)+1,\  2\le i<k/2, \ 2S(k/2)-1).$$
\end{prop}

To be consistent with property (3) of \fullref{gthm}, our function $G$ must satisfy the property stated in the next theorem,
the proof of which will occupy much of this section.

\begin{thm}\label{firstproperty} If $e\ge8$, and $2\le i,j\le2^{e-4}$, then
\begin{equation}G(e-1,i)+G(e-1,j)+1-G(e,i+j)\ge0.\label{Gineq}\end{equation}
\end{thm}

The stabilization given by $S(k)$ and the requirement \eqref{Gineq} are what lead to the complicated sum in \eqref{geq1}.
The first example of this is for $G(11,3)+G(11,32)+1-G(12,35)$. Since $G'(11,3)=9>S(3)=8$, we have $G(11,3)=8$. Also $G(11,32)=11-6-5=0$,
and $\la 12-6-\nu'(35)\ra=10$. Thus we must subtract 1 from $\la 12-6-\nu'(35)\ra$ in $G(12,35)$
in order that \eqref{Gineq} will hold. This is accounted for by
the decomposition of $35$ with $k_1=3$. The value $R(3)=4$ is the amount that $\nu(k-3)$ must exceed in order that the decomposition
affects the value of $G(e,k)$.

Note that 11 is the smallest value of $e$ for which $G(e,3)\ne\la e-6-\nu'(3)\ra$. This is obtained by solving
$$e-6-\nu'(3)=S(3)+1,$$
obtaining $e=11$.
We want $R(3)$ to be 1 less than the value of $t$ which satisfies
$$G(11,3)+G(11,2^t)+1-\la 12-6-\nu'(2^t+3)\ra=-1.$$
Here $G(11,3)-\la 12-6-\nu'(2^t+3)\ra$ necessarily equals $-2$:\ 1 from $12-11$, and 1 from $G(11,3)=G'(11,3)-1$. Thus we need $t$ to satisfy
$0=G(11,2^t)=11-6-t$, and so
$$R(3)=t-1=(S(3)+\nu'(3)+6)-6=S(3)+\nu'(3),$$
consistent with \eqref{Rdef}, since $V(3)=0$.

The way $V$ arises can be seen by comparing the requirements, for $t\ge5$,
$$G(e,2^t+5)\le G(e-1,2)+G(e-1,2^t+3)+1$$
and
$$G(e,2^t+5)\le G(e-1,5)+G(e-1,2^t)+1.$$
The first reduces to, for $e$ moderately large,
$$G(e,2^t+5)\le S(2)+e-6-t+S(3)=e+6-t,$$
while the second becomes
$$G(e,2^t+5)\le S(5)+e-6-t=e+7-t.$$
We must use the first condition because $S(2)+S(3)<S(5)$. The value $V(5)=1$ measures this. Our $V(k)$ satisfies that it is the largest
$r$ such that $k=i_0+\cdots+i_r$ with
$$S(i_0+\cdots+i_t)=S(i_0+\cdots+i_{t-1})+S(i_t)+1$$
for $1\le t\le r$.

This concludes our discussion of the rationale behind the definition of $G$ except for one more brief comment.
It was certainly to be expected that these modifications to the $G$-formula, given by the summands in \eqref{geq1},
would be cumulative. It was not {\it a priori} clear whether $R(k_{i-1})$ or $R(k_0+\cdots+k_{i-1})$ would be the appropriate part
of that formula. The answer will become apparent in Subcase 2d of the proof of \fullref{not8k+1}.

The following proposition will be needed shortly. The function $S'$ below will often be encountered in the
guise of $S'(k)=R(k)-\nu'(k)$.
\begin{prop} Let $S'(k)=S(k)-V(k)$. If $i,j,i+j\not\equiv1$ mod $8$, then
$$S'(i)+S'(j)\ge S'(i+j).$$\label{spprop}
Moreover, if $i<2^{\nu(j)}$, then equality is obtained.
\end{prop}
\begin{proof}
One easily verifies that
$$S'(k)=8k-13[\tfrac{k+1}2]+\a(k)+\begin{cases}0&k\equiv0\ (2)\\
8&k\equiv3\ (4)\\
9&k\equiv5\ (8).\end{cases}$$
For $1\le m\le 4$, let $\phi_m$ denote the $m$th part of the above formula for $S'(k)$, and let $\psi_m(i,j)=\phi_m(i)+\phi_m(j)-\phi_m(i+j)$.
Then $$\psi_m(i,j)=\begin{cases}0&m=1\\
0&m=2,\ ij\text{ even}\\
-13&m=2,\ ij\text{ odd}\\
\nu\tbinom{i+j}i&m=3\\
\ge16&m=4,\ ij\text{ odd}\\
-1&m=4,\ i+j\equiv5 \ (8)\text{ and }i\text{ or }j\equiv 3\ (4)\\
\ge0&m=4, \text{otherwise.}\end{cases}$$
Since $\binom{i+j}{i}$ is even if $i+j\equiv5\ (8)$ and $i\equiv3\ (4)$, the inequality follows.

For the second part, one easily sees that, if $i<2^{\nu(j)}$, then $\psi_m(i,j)=0$ for $1\le m\le4$. When $m=4$, it is true
because $i\equiv i+j$ mod 8 (or $\nu(j)=2$ and $i=2$ or 3).
\end{proof}

We now begin the lengthy proof of \fullref{firstproperty}. In order to keep the number of cases and subcases within
reason, we split the theorem into two parts. Most of the work will go into proving the following result.
\begin{thm}\label{not8k+1} If $e\ge8$,  $2\le i,j\le2^{e-4}$, and $i,j,i+j\not\equiv1$ mod 8, then
\eqref{Gineq} holds.
\end{thm}
\begin{proof} We divide into cases depending upon whether $S(i)$ and/or decompositions are involved.
\ss

\noindent {\bf Case 1}: Neither $i$ nor $j$ decomposes,  $G(e-1,i)\ne S(i)$, and $G(e-1,j)\ne S(j)$.
In this case, the LHS of \eqref{Gineq} becomes
\begin{equation}\label{lara}\ge\la e-7-\nu'(i)\ra+\la e-7-\nu'(j)\ra+1-\la e-6-\nu'(i+j)\ra.\end{equation}
By considering separately the four subcases (a) $i$ and $j$ odd, (b) $i$ odd, $j$ even, (c) $\nu(j)>\nu(i)>0$, and
(d) $\nu(i)=\nu(j)>0$, one easily shows that \eqref{lara} is $\ge0$ in each subcase. Note that if $i+j$ decomposes,
then the LHS of \eqref{Gineq} is greater than \eqref{lara}, and so we need not worry about this possibility here.

\ss
\noindent {\bf Case 2}: $G(e-1,i)=S(i)$ and $i$ does not decompose.

{\bf Subcase 2a}: Also, $G(e-1,j)=S(j)$. Then the LHS of \eqref{Gineq} is $\ge S(i)+S(j)+1-S(i+j)\ge0$, by \fullref{Sprop}. The remaining subcases of Case 2 now assume that $G(e-1,j)<S(j)$.

\ss
{\bf Subcase 2b}: $j$ does not decompose, and $\nu(j)\le\nu(i)$. Then $\nu'(i+j)\ge\nu'(j)$, and so the LHS of \eqref{Gineq} is
$$\ge S(i)+\la e-7-\nu'(j)\ra+1-\la e-6-\nu'(i+j)\ra\ge S(i)>0.$$

\ss
{\bf Subcase 2c}: $j$ does not decompose, and $\nu(j)>\nu(i)$. We allow for the possibility that $i$ might
serve as the bottom part of a decomposition of $i+j$. This will be true if $\nu(j)$ is sufficiently large.
Because of our $\la-\ra$-notation, our analysis is valid regardless.
This time $\nu'(i)=\nu'(i+j)$, and so the LHS of \eqref{Gineq} is
$$\ge S(i)+\la e-7-\nu(j)\ra+1-\la e-6-\nu'(i)\ra+\la\min(\nu(j),e-6)-S(i)-\nu'(i)+V(i)\ra.$$
If $\nu(j)\le e-7$, this is $\ge V(i)\ge0$. If $\nu(j)\ge e-6$, it simplifies to
\begin{equation}\label{2c}\ge V(i)+1+e-6-\nu'(i)-\la e-6-\nu'(i)\ra.\end{equation}
Since $j\le 2^{e-4}$ and $\nu(i)<\nu(j)$, we have $\nu'(i)\le e-5$, and so \eqref{2c} is $\ge V(i)\ge0$.

\ss
{\bf Subcase 2d}: $j$ admits a decomposition. We  consider a 2-stage decomposition $j=j_0+j_1+2^tA$ with
$A$ odd and $\nu(j_1)>R(j_0)$.
It will be clear that the argument here can be adapted to a longer decomposition. Letting $D\ge0$ denote
any amount added for a decomposition of $i+j$,
the LHS of \eqref{Gineq} becomes, using \fullref{spprop},
\begin{eqnarray}&&S(i)+( e-7-\nu'(j))-(\nu(j_1)-R(j_0))\nonumber\\
&&\quad-(\min(t,e-7)-R(j_1))+1-\la e-6-\nu'(i+j)\ra+D\label{minD}\\
&=&S'(i)+V(i)+S'(j_0)+S'(j_1)+\nu'(i+j)-\min(t,e-7)+D\nonumber\\
&\ge&V(i)+S'(i+j_0+j_1)+\nu'(i+j)-\min(t,e-7)+D\nonumber\\
&=&V(i)+R(i+j_0+j_1)-\nu'(i+j_0+j_1)+\nu'(i+j)-\min(t,e-7)+D\label{2d}.\end{eqnarray}
We will discuss later the removal of the $\la-\ra$ at the first step.

We will show below that \begin{equation}\label{V}V(i)-\nu'(i+j_0+j_1)+\nu'(i+j)>0.\end{equation}
Assuming this, the only way that \eqref{2d} could be negative is if $\min(t,e-7)>R(i+j_0+j_1)$.
But if this is the case, then $(i+j_0+j_1)+2^tA$ is a decomposition of $i+j$, which makes
 $D\ge\min(t,e-6)-R(i+j_0+j_1)$. If $i+j_0+j_1$ decomposes
further, that only adds more to $D$. Thus, assuming \eqref{V}, we obtain that \eqref{2d} is $\ge0$.

We now prove \eqref{V}. The only way it could possibly be negative is if $i=2^tB-j_0-j_1$ with
$B$ even. Then the LHS of \eqref{V} becomes
\begin{eqnarray*}&\ge& \a(2^tB-j_0-j_1)-2-(t+\nu(B))+t\\
&=&\a(B-1)+t-\a(j_0+j_1-1)-2-\nu(B)\\
&>&0\end{eqnarray*}
since $\a(B-1)\ge\nu(B)$ and $t\ge R(j_0+j_1)>>\a(j_0+j_1-1)$.

Regarding the removal of $\la-\ra$ above: if $\nu'(i+j)>e-6$, then  \eqref{minD} becomes
\begin{eqnarray*}&\ge&S(i)+e-7-\nu'(j)-\nu(j_1)+R(j_0)-\min(t,e-7)+R(j_1)+1\\
&=&S(i)+(e-7-\min(t,e-7))+(R(j_1)-\nu(j_1))+(R(j_0)-\nu'(j_0))+1\\
&>&0\end{eqnarray*}
because each of its terms is nonnegative.

\ss
\noindent{\bf Case 3}: $G(e-1,i)=S(i)$ and $i$ decomposes. Although the decomposition
of $i$ does not affect the value of $S(i)$, it could affect the value of $G(e,i+j)$
by affecting the decomposition of $i+j$. In the analogues of Subcases 2a and 2b, the
decomposition of $i+j$ was not needed, and so a decomposition of $i$ cannot affect the validity.

\ss
{\bf Subcase 3a}: $j$ does not decompose and $\nu(j)>\nu(i)$.

\ss
{\bf Subsubcase 3ai}: $i$ admits a simple decomposition. Let $i=i_0+2^t\a$ with $\a$ odd
and $t>R(i_0)$.
If $\nu(j)\ge t$, then, considering $i_0+(2^t\a)+j$ as a possible decomposition of $i+j$, the LHS of \eqref{Gineq} becomes
\begin{eqnarray*}&\ge&S(i)+\la e-7-\nu(j)\ra+1-\la e-6-\nu'(i)\ra\\
&&\quad +\la\min(t,e-6)-R(i_0)\ra+\la\min(\nu(j),e-6)-R(2^t\a)\ra.\end{eqnarray*}
This exceeds the amount analyzed in Subcase 2c by
\begin{equation}t-R(i_0)-R(2^t\a)+R(i_0+2^t\a).\label{tR}\end{equation}
Since, in the notation of \fullref{spprop}, $S'=R-\nu'$, and $\nu'(i_0)=\nu'(i_0+2^t\a)$, then
\eqref{tR} equals $S'(i_0+2^t\a)-S'(i_0)-S'(2^t\a)=0$ by \fullref{spprop}.

If, on the other hand, $\nu(j)<t$, then we don't need $i+j$ to be decomposable, since the LHS of
\eqref{Gineq} $$\ge S(i)+\la e-7-\nu(j)\ra+1-\la e-6-\nu'(i)\ra=S(i)+\nu'(i)-\nu(j)>0,$$
since $S(i)>t+4>\nu(j)+4$. (The $+4$ is included because of the possibility that $\nu'(i)=-4$.)

\ss
{\bf Subsubcase 3aii}: $i$ admits a multiple decomposition. If $\nu(j)\le S(i)+\nu'(i)$, then, as in the preceding paragraph,
we do not need a decomposition of $i+j$ in order to satisfy \eqref{Gineq}. If, on the other hand, $\nu(j)>S(i)+\nu'(i)$, then
the result follows as in the first paragraph of Subcase 3ai, using additivity of $S'$ on disjoint decompositions.

\ss
{\bf Subcase 3b}: $i$ and $j$ both decompose exactly once. Let $i=i_0+2^m\b$ with $\b$ odd and $m> R(i_0)$, and $j=j_0+2^t\a$ with $\a$ odd and
$t> R(j_0)$.

If $m>t$, then we can consider $i+j$ as
$(i_0+j_0)+2^t\a+2^m\b$. It is possible that $\la m-R(2^t\a)\ra$ might contribute to $G(e,i+j)$, but even if it does, we do not need it.
The situation is similar to Subcase 2d. Using the $\la t-R(j_0)\ra$ and $\la t-R(i_0+j_0)\ra$ parts of $G(e-1,j)$ and $G(e,i+j)$,
respectively, the LHS of \eqref{Gineq} simplifies to
$$\ge S(i)-\nu'(j_0)+\nu'(i_0+j_0)+R(j_0)-R(i_0+j_0),$$
which is very positive. (It would be $\ge V(i_0)$ by \fullref{spprop} if $S(i)$ were replaced by the much smaller number $S(i_0)$.)
Keeping in mind that $2^{e-3}\ge i+j$, we will usually omit, from now on, explicit consideration of the possibility that $e-6<\nu(k_{i-1})$ in \eqref{geq1}.
In Subcase 4d, there is a detailed discussion of a delicate case in which we consider carefully what happens when $e-6$ is larger than
the relevant 2-exponent .

If $m=t$, then a very similar argument works. Because the decomposition of $i+j$ now is $(i_0+j_0)+2^p\gamma$ with $p>t$,
and this exponent appears with a $+$ sign in $-G(e,i+j)$, the LHS of \eqref{Gineq} is even larger than it was when $m>t$.

Now suppose $m<t$. We use $(i_0+j_0+2^m\b)+2^t\a$ as our trial decomposition of $i+j$. If it is not a true
decomposition, then the $\la-\ra$ will take care of it.

The LHS of \eqref{Gineq} becomes
\begin{eqnarray*} &\ge&S(i)+(e-7-\nu'(j_0)) -(t-R(j_0))+1\\
&&\ -(e-6-\nu'(i_0+j_0))+\la t-R(i_0+j_0+2^m\b)\ra\\
&\ge&S(i)-\nu'(j_0)+R(j_0)+\nu'(i_0+j_0)-R(i+j_0)\\
&=&V(i)+S'(i)+S'(j_0)-S'(i+j_0)\\
&\ge&V(i).\end{eqnarray*}


{\bf Subcase 3c}: At least one of $i$ and $j$ decomposes more than once. The argument is very similar to that of Subcase 3b.
The only reason for separating them is to use 3b as a warmup for 3c.
Let $i=i_0+\cdots+i_r$ and $j=j_0+\cdots+j_s$ be maximal decompositions.

If $\nu(j_s)\le \nu(i_r)$, then the LHS of \eqref{Gineq} is, without using any decomposition of $i+j$,
\begin{eqnarray*} &\ge&S(i)-\nu'(j)-\sum_{k=1}^s(\nu(j_k)-R(j_{k-1}))+\nu'(i+j)\\
&\ge&S(i)+\sum_{k=0}^{s-1}(R(j_k)-\nu'(j_k))-\nu(j_s)\\
&\ge& S(i)-\nu(i_r)\\
&>\!>&0.\end{eqnarray*}

If $\nu(i_r)<\nu(j_s)$, first suppose the only decomposition of $i+j$ is the simple decomposition $K+j_s$ with $K=i+j_0+\cdots+j_{s-1}$.
Then the LHS of \eqref{Gineq} is
\begin{eqnarray*}&\ge&S(i)-\nu'(j)-\sum_{k=1}^s(\nu(j_k)-R(j_{k-1}))+\nu'(i+j)+\nu(j_s)-R(K)\\
&=&R(i)+V(i)-\nu'(i)+\sum_{k=0}^{s-1}(R(j_k)-\nu'(j_k))+\nu'(K)-R(K)\\
&\ge&V(i)\end{eqnarray*}
by \fullref{spprop}.

If $i+j$ decomposes more finely, say as $A+B+j_s$, then $-R(K)$ is replaced by $-R(B)+\nu(B)-R(A)$. But these are equal by the
second part of \fullref{spprop}, noting that $\nu'(A+B)=\nu'(A)$.

\medskip
{\bf Case 4}: $S(-)$ not involved, $i$ decomposes, $j$ doesn't. Recall $i,j\le 2^{e-4}$. We assume that $i$ admits a decomposition
as $i_0+i_1+i_2$. The nature of our argument will show that the conclusion will also be true for longer
decompositions. The LHS of \eqref{Gineq} becomes
\begin{equation} \label{case4} e-6-\nu'(i_0)-\nu(i_1)+R(i_0)-\nu(i_2)+R(i_1)+\la e-7-\nu'(j)\ra+Y,\end{equation}
where $Y=-G(e,i+j)$. We use \eqref{X} often in what follows.

\ss
{\bf Subcase 4a}: $\nu(j)<\nu(i)$. Then, using a decomposition $i+j=(i_0+j+i_1)+(i_2)$, we obtain
\begin{equation}\label{-1}Y\ge-(e-6-\nu'(j))+\la\nu(i_2)-R(i_0+j+i_1)\ra.\end{equation}
If there is an additional decomposition of $i+j$ as $(i_0+j)+(i_1)+(i_2)$, then by the second part of
\fullref{spprop}, $R(i_0+j+i_1)=R(i_0+j)+R(i_1)-\nu(i_1)$, and so the same expression is obtained.
Then \eqref{case4} is
\begin{equation}\label{-2}\ge(e-7-R(i_0+j+i_1))+(R(i_0)-\nu'(i_0))+(R(i_1)-\nu(i_1))>0,\end{equation}
since if the $\la-\ra$ in \eqref{-1} is $>0$, then $$e-7\ge\nu(i_2)-2\ge R(i_0+j+i_1)-2,$$
but the $(R-\nu)$-expressions are $>2$. If the $\la-\ra$ in \eqref{-1} is 0, then the first
part of \eqref{-2} is replaced by $(e-7-\nu(i_2))\ge-2$.

\ss
{\bf Subcase 4b}: $\nu(i)\le \nu(j)<R(i_0)$. In this case, which is very similar to 4a,
$$Y\ge-(e-6-\nu'(i))+\la\nu(i_2)-R(i_0+j+i_1)\ra,$$
because if there is an additional decomposition of $i+j$ as $(i_0+j)+(i_1)+(i_2)$, then
$R(i_0+j+i_1)=R(i_0+j)+R(i_1)-\nu(i_1)$, and so the expression for $Y$ is unchanged.
Then \eqref{case4} is
$$\ge(S'(i_0)+S'(i_1)+S'(j)-S'(i_0+i_1+j))+(e-7-R(j))>0.$$
In the remaining subcases, we deal with a maximum possible decomposition of $i+j$, realizing, as in 4a and 4b,
that if the decomposition must be amalgamated, the expression is not changed.

\ss
{\bf Subcase 4c}: $R(i_0)\le \nu(j)<\nu(i_1)$. Then
$$Y\ge-(e-6-\nu'(i))+\la\nu(j)-R(i_0)\ra+\la\nu(i_1)-R(j)\ra+\la\nu(i_2)-R(i_1)\ra,$$
and so \eqref{case4} is
$$\ge(e-7-\nu(i_1))+\la\nu(i_1)-R(j)\ra>0.$$

\ss
{\bf Subcase 4d}: $\nu(i_1)\le \nu(j)<\nu(i_2)$. Then
$$Y\ge-(e-6-\nu'(i))+\la\nu(i_1)-R(i_0)\ra+\la\nu(j)-R(i_1)\ra+\la\nu(i_2)-R(j)\ra,$$
and so \eqref{case4} is
\begin{equation}\label{4d}\ge(e-7-\nu(i_2))+\la\nu(i_2)-R(j)\ra>0.\end{equation}

As noted in Subcase 3b, we are usually not paying explicit attention to the possibility
that $e-6\le\nu(i_2)$ (in the situation in this subcase, 4d). Here it does warrant our
attention. We might have $i_2=2^{e-5}$, $2^{e-6}$, or $3\cdot2^{e-6}$, and then it would
seem that \eqref{4d} might not be valid.

If $i_2=2^{e-5}$, then $\la \nu(i_2)-R(i_1)\ra$ in the above analysis is replaced by
$\la e-7-R(i_1)\ra$. This decrease of 2 compensates for the fact that $e-7-\nu(i_2)=-2$
in \eqref{4d}. Similarly, if $\nu(i_2)=e-6$, then $\la\nu(i_2)-R(i_1)\ra$ is replaced by
$\la e-7-R(i_1)\ra$, compensating for $e-7-\nu(i_2)=-1$.

\ss
{\bf Subcase 4e}: $ \nu(i_2)<\nu(j)$. Then
$$Y\ge-(e-6-\nu'(i))+\la\nu(i_1)-R(i_0)\ra+\la\nu(i_2)-R(i_1)\ra+\la\nu(j)-R(i_2)\ra,$$
and so \eqref{case4} is
$$\ge\la e-7-\nu(j)\ra+\la\nu(j)-R(i_2)\ra>0.$$

\ss
{\bf Case 5}: $S(-)$ not involved, both $i$ and $j$ decompose.  We consider here a typical
example in which both $i$ and $j$ decompose twice. It should be clear that the general case
will work out in the same way. We assume that $i=i_0+i_1+i_2$ and $j=j_0+j_1+j_2$
are decompositions. Then
\begin{eqnarray*}G(e-1,i)+G(e-1,j)+1&=&
e-6-\nu'(i_0)-\nu(i_1)+R(i_0)-\nu(i_2)+R(i_1)\\&&+e-7-\nu'(j_0)-\nu(j_1)
+R(j_0)-\nu(j_2)+R(j_1)
\end{eqnarray*}
We assume without much loss of generality that $\nu(j_2)>\nu(i_2)$ and $\nu(i_0)<\nu(j_0)$.

\ss
{\bf Subcase 5a}: $\nu(j_2)<R(i_0+i_1+i_2+j_0+j_1)$. We use no decomposition of $i+j$.
We obtain that
\begin{eqnarray*}&&G(e-1,i)+G(e-1,j)+1-G(e,i+j)\\
&\ge&R(i_0)+S'(i_1)-\nu(i_2)+e-7+S'(j_0)+S'(j_1)-\nu(j_2)\\
&=&S'(i_0)+S'(i_1)+S'(i_2)+S'(j_0)+S'(j_1)+\nu'(i_0)-R(i_2)+e-7-\nu(j_2)\\
&\ge&R(i_0+i_1+i_2+j_0+j_1)-R(i_2)+e-7-\nu(j_2)\\
&>\!>&0,\end{eqnarray*}
since $e-7-\nu(j_2)\ge-2$ while $R(i_0+i_1+i_2+j_0+j_1)-R(i_2)>\!>0$.

\ss
{\bf Subcase 5b}: $\nu(j_2)>R(i_0+i_1+i_2+j_0+j_1)$. We use a decomposition of $i+j$
as $(i_0+i_1+i_2+j_0+j_1)+(j_2)$. We discuss afterward the usual argument regarding what
happens if it decomposes more finely. Similarly to Subcase 5a, we obtain
\begin{eqnarray*}&&G(e-1,i)+G(e-1,j)+1-G(e,i+j)\\
&\ge&S'(i_0)+S'(i_1)+S'(i_2)+S'(j_0)+S'(j_1)-S'(i_0+i_1+i_2+j_0+j_1)\\
&&\ +e-7-R(i_2)\\
&\ge&0\end{eqnarray*}
using \fullref{spprop} and $$e-7\ge\nu(j_2)-2\ge R(i_0+i_1+i_2+j_0+j_1)-2>\!>R(i_2).$$

Further decomposition of $i_0+i_1+i_2+j_0+j_1$ into 2-adically disjoint parts does not
change the expression, using the second part of \fullref{spprop}, similarly to the argument
in Subcases 4a and 4b.
\end{proof}

The following result will be useful in some subsequent proofs. In particular,
\fullref{gthm}.(5) is an immediate consequence.
\begin{prop}\label{minus1} For $e\ge7$ and $2\le k<2^{e-3}$,
$$G(e,k+1)-G(e,k)\begin{cases}=8&k\equiv0\ (8),\ \a(k)=1,\ e\ge S(k)+\nu(k)+8\\
=7&k\equiv0\ (8),\ \a(k)=1,\ e= S(k)+\nu(k)+7\\
=6&k\equiv0\ (8),\text{ otherwise}\\
\le-1&k\equiv1\ (8)\\
\le6&\text{otherwise.}\end{cases}$$
\end{prop}
\begin{proof} We begin by noting that the result is true for the limiting values, $S(k)$,
since they are easily shown to satisfy
\begin{equation}S(k+1)-S(k)\begin{cases}=8&k=2^e,\ e\ge3\\
=6&k\equiv0\ (8),\ \a(k)>1\\
=6&k\equiv 4\ (8)\\
=4&k\equiv2\ (4)\\
\le-1&k\equiv1\ (8)\\
=-1&k\equiv3,5\ (8)\\
\le-3&k\equiv7\ (8).\end{cases}\label{Scases}\end{equation}
The case $k\equiv0$ mod 8 of the proposition follows easily from \eqref{Scases}
and the definitions.

We next handle the case $k=8\ell+1$. If $\nu(\ell)\ge3$, then $8\ell+2$ admits a
decomposition with $k_0=2$. Any additional portions of a decomposition of $8\ell+2$
will occur identically in $8\ell$. Thus, in this case, with $\nu=\nu(8\ell)\ge6$,
$$G(e,8\ell+2)-G(e,8\ell+1)=e-7-\la\min(\nu,e-6)-5\ra-(6+\la e-6-\nu\ra).$$
This is $\le-2$, regardless of the sign of $e-6-\nu$.

Now assume $\nu(\ell)<3$. If $8\ell$ admits a decomposition as $k_0+2^t\a$ with $\a$ odd,
then we consider $(k_0+2)+2^t\a$ as a possible decomposition of $8\ell+2$. Any additional
portions of a decomposition of $8\ell+2$ occur identically in $8\ell$. For $v=\nu(\ell)=0$, 1, or 2,
we obtain
\begin{equation}G(e,8\ell+2)-G(e,8\ell+1)=e-13-\la e-9-v\ra-\la D-2+v\ra+\la D\ra,\label{DE}\end{equation}
where $D=\min(t,e-6)-R(k_0)$. Here we have used the easily-verified fact that
if $k_0\equiv0$ mod 8, then $R(k_0+2)-R(k_0)=5-\nu(k_0)$. One easily checks that
\eqref{DE} is $\le -2$ for any $e$ and $D$, since $0\le v\le 2$.

For $\tau=[2,3,4,5,6,7]$ and $k=8\ell+\tau$, we have, for $e>7$,
$$\la e-6-\nu'(k+1)\ra-\la e-6-\nu'(k)\ra=[5,-6,6,-5,5,\le-5],$$
and, if $k$ admits a simple decomposition $k_0+2^t\a$ with $\a$ odd,
$$\la m-R(k_0)\ra -\la m-R(k_0+1)\ra\le[0,5,0,4,0,5].$$
Here $m=\min(e-6,t)$. As before, higher deviations will cancel in the difference.
Thus $G(e,k+1)-G(e,k)$, which is the sum of the two displays of this paragraph, is $\le6$,
as claimed.
\end{proof}

Now we can complete the proof of \fullref{firstproperty} by proving.
\begin{thm}\label{1mod8} Theorem \ref{firstproperty} is true when
$i$ or $j$ or $i+j$ is $\equiv1$ mod $8$.\end{thm}
\begin{proof} Again we divide into cases.

\ss
{\bf Case 1}: Only $i\equiv1$ mod 8. We have
\begin{eqnarray}&&\nonumber G(e-1,i)+G(e-1,j)+1-G(e,i+j)\\
&=&\bigl(G(e-1,i)-G(e-1,i-1)\bigr)-\bigl(G(e,i+j)-G(e,i-1+j)\bigr)\nonumber\\
&&+\bigl(G(e-1,i-1)+G(e-1,j)+1-G(e,i-1+j)\bigr)\label{diff}\\
&\ge&0,\nonumber\end{eqnarray}
since the first $(-)$ in \eqref{diff} is $\ge6$ by \fullref{minus1}, the second is
$\le6$ by \fullref{minus1}, and the third is $\ge0$ by \fullref{not8k+1}.

\ss
{\bf Case 2}: both $i$ and $j \equiv1$ mod 8. This follows by an argument similar to
that of Case 1.

\ss
{\bf Case 3}: $i$ and $i+j\equiv1$ mod 8. This follows from the validity for $(i-1,j)$
similarly to Case 1. Usually $G(e-1,i)-G(e-1,i-1)=6$ and $G(e,i+j)-G(e,i-1+j)=6$, and so the inequality follows as in \eqref{diff}. If $G(e,i+j)-G(e,i-1+j)>6$, then $G(e-1,i)=S(i)$ and $G(e-1,j)=S(j)$, and so
$$G(e-1,i)+G(e-1,j)+1-G(e,i+j)\ge S(i)+S(j)+1-S(i+j)\ge0$$
by \fullref{Sprop}.

\ss
{\bf Case 4}: $i+j\equiv1$ mod 8, while $i,j\not\equiv1$ mod 8. If $G(e,i+j)-G(e,i+j-1)>6$, then
$G(e-1,i)=S(i)$, $G(e-1,j)=S(j)$, and $G(e,i+j)\le S(i+j)$, and so the result follows from \fullref{Sprop}.
So we may now assume $G(e,i+j)-G(e,i+j-1)=6$.
Without loss of generality,
assume $i$ is odd and $j$ is even.

First, we assume $i\equiv3$ mod 4. By the proof of \fullref{minus1},
$G(e,i)-G(e,i-1)=4$ or 5, and if $i$ is indecomposable, then $G(e,i)-G(e,i-1)=
4$ if and only if $G(e,i-1)=S(i-1)$.
Thus the result will follow as in \eqref{diff} once we show that if $i,j\equiv2$ mod 4 and $i+j\equiv0$ mod 8, then \eqref{Gineq} is satisfied with
1 to spare, and with 2 to spare if $G(e,i)-G(e,i-1)=4$.

The basic value of the LHS of \eqref{Gineq} in this case is
\begin{equation}\label{e-8}\la e-8\ra+\la e-8\ra+1-\la e-v\ra\end{equation}
with $v\ge9$. This equals 1 if $e=7$ or 8, while for $e\ge9$, it is $\ge e-6$.
The smallest $e$ for which the LHS of \eqref{Gineq} does not equal \eqref{e-8}
is $e=12$, when $i=2$.

Neglecting temporarily the effect of deviations, the desired conclusion is obtained
since it is true at the onset of $S(i)$ and will continue to be true as $e$ increases,
since now $G(e-1,j)$ and $G(e,i+j)$ will both increase by 1 each time. When $G(e-1,j)$ achieves
a value of $S(j)$, then the LHS of \eqref{Gineq} is
$$\ge S(i)+S(j)+1-S(i+j)>2$$
for the congruences being considered here.

When deviations are taken into account, the fact that makes it work is the easily-verified fact that
\begin{equation}\label{26}R(8\ell+2)+R(8\ell'+6)-R(8\ell+8\ell'+8)=1+\nu\bigl(\tbinom{\ell+\ell'}{\ell}\bigr).
\end{equation}
Suppose, for example, that $i=i_0+2^t\a$ and $j=j_0+2^u\b$ are decompositions with $\a$ and $\b$ odd, and $t<u\le e-7$.
The LHS of \eqref{Gineq} becomes
$$\ge e-8+e-7-(t-R(i_0))-(u-R(j_0))-(e-v)+\la t-R(i_0+j_0)\ra+\la  u-R(2^t\a)\ra$$
with $v\ge9$. Using \eqref{26}, this is
\begin{eqnarray*}&\ge&e+v-14+R(i_0+j_0)-t+\nu\bigl(\tbinom{i_0+j_0}{i_0}\bigr)+\la t-R(i_0+j_0)\ra-u\\
&\ge& v-7\\
&\ge&2,\end{eqnarray*}
since $e-7\ge u$ and using \eqref{X}.
Other situations involving decompositions work out similarly.

The case $i\equiv5$ is handled similarly.
\end{proof}

Next we verify the first part of \fullref{gthm}.(4). In fact the conclusion of that theorem is true without
regard for the hypothesis.
\begin{thm}\label{nohyp} If $i,j\le 2^{e-3}$ and $i+j+1\le2^{e-2}$, then
$$g(e,i)+g(e,j)<g(e+1,i+j+1).$$
\end{thm}
\begin{proof} We prove the equivalent statement, with $i$, $j$, and $e$ as in the hypothesis,
\begin{equation}\label{equiv}G(e,i)+G(e,j)+8>G(e+1,i+j+1).\end{equation}
By \fullref{firstproperty} and \fullref{minus1}, we have
$$G(e,i)+G(e,j)+8\ge G(e+1,i+j)+7>G(e+1,i+j+1)$$
unless $i+j+1=2^t+1$ with $t\ge3$ and $G(e+2,i+j+1)=S(i+j+1)$. In this case, it will
also be true that $G(e,i)=S(i)$ and $G(e,j)=S(j)$. Thus it suffices to show
$$S(i)+S(2^t-i)+8>S(2^t+1).$$
This follows readily from the definition of $S$. The smallest value of $S(i)+S(2^t-i)$ occurs
when $i=2^{t-1}$ and is $3\cdot2^{t-1}+2$, while $S(2^t+1)=3\cdot2^{t-1}+9$.
\end{proof}

The second part of \fullref{gthm}.(4) follows from the following result.
\begin{thm}\label{2k} For $k\le 2^{e-3}$, $G(e+1,2k)\le 2G(e,k)$ with equality if and only if
$G(e+1,2k)=G(e,k)=0$, which occurs if and only if $$k\in\{2^{e-3},\ 2^{e-4}, \ 2^{e-5},\ 3\cdot2^{e-5},\ 2^{e-6}\a\}\text{ with }\a\in\{1,3,5,7\}.$$
If equality occurs, then
$$G(e+1,2k)<G(e,\ell)+G(e,2k-\ell)+1$$
for all $\ell$.
\end{thm}
\begin{proof} The second sentence follows immediately from the first, since
$$0<G(e,\ell)+G(e,2k-\ell)+1.$$

For basic values, we have
$$2G(e,k)-G(e+1,2k)=\begin{cases}2(e-2)-(e-6)&k\text{ odd}\\
\la e-6-\nu(k)\ra&k\text{ even.}\end{cases}$$
This is clearly $\ge0$, and $=0$ in exactly the cases claimed.

If $G(e,k)=S(k)$, then
$$2G(e,k)-G(e+1,2k)\ge 2S(k)-S(2k)=\begin{cases}2\a(k)-1&k\text{ even}\\
12&k=2^t+1,\ t\ge3\\
2\a(k)+4&k\equiv1\ (8),\ \a(k)\ne2\\
4\nu(k-1)+2\a(k)-4&k\equiv3,5,7\ (8).\end{cases}$$
This is $>0$.

Suppose $k=k_0+2^t\a$ is a simple decomposition, with $\a$ odd and $e-6\ge t$.
If $k$ is even, then $2R(k)=R(2k)+\a(k-1)$, and so
\begin{eqnarray*}&&2G(e,k)-G(e+1,2k)\\
&=&e-6-\nu(k)-2(t-R(k_0))+\la t+1-R(2k_0)\ra\\
&=&e-5-t+\a(k_0)-1+R(2k_0)-t-1+\la t+1-R(2k_0)\ra\\
&\ge&1,\end{eqnarray*}
using \eqref{X}.
If $k\equiv3,5,7$ mod 8, then
$$2R(k)=R(2k)+4\nu(k-1)+\a(k)-\begin{cases}10&k\equiv3\ (4)\\ 12&k\equiv1\ (4).\end{cases}$$
Then
\begin{eqnarray*}&&2G(e,k)-G(e+1,2k)\\
&=&e+2-2(t-R(k_0))+\la t+1-R(2k_0)\ra\\
&\ge& e+3-t+R(2k_0)-t-1+\la t+1-R(2k_0)\ra+4\nu(k_0-1)+\a(k_0)-12\\
&>&0.\end{eqnarray*}
The situation when $t>e-6$ and the case of higher deviations are handled similarly.

Finally, we have
\begin{eqnarray*}&&2G(e,8\ell+1)-G(e+1,16\ell+2)\\
&\ge&2(G(e,8\ell)+6)-G(e+1,16\ell)-(G(e+1,16\ell+1)-G(e+1,16\ell))\\
&&\ -(G(e+1,16\ell+2)-G(e+1,16\ell+1))\\
&\ge&12+0-8-(-1)\\
&>&0.\end{eqnarray*}
\end{proof}

Finally, we verify part (2) of \fullref{gthm}. We have
$$g(e,k)=8k-G(e,k)\ge 8k-S(k)\ge13[\tfrac{k+1}2]-2\a(k)-10.$$
This is $\ge 4k+4$ for $k\ge7$, while for $k<7$ we verify directly that $8k-S(k)\ge 4k+4$.

\section{A bound for geometric dimension of normal bundle}\label{normsec}
In this section, we prove the following key result, a main ingredient in the proof of our geometric dimension result, \fullref{1better}, which has
already been seen to imply our immersion theorem.
\begin{thm} If $e\ge7$ and $t\ge 1$ and $k_0+\cdots+k_{t-1}=2^{e-3}-1$, then \begin{equation}\label{immsum}\sum_{i=0}^{t-1} G(e+i,k_i)\ge e-2.
\end{equation}\label{gdsum}\end{thm}

\begin{rmk} {\rm The integers $k_i$ in this theorem are nonnegative, but possibly zero. Some examples in which equality is obtained are}
\begin{itemize}
\item $G(e,2^{e-3}-1)$;
\item $G(e,2^{e-4}-1)+G(e+1,2^{e-4})$;
\item $G(e,2^{e-5}-1)+G(e+1,3\cdot2^{e-5})$;
\item $G(e,3\cdot2^{e-5}-1)+G(e+1,2^{e-5})$;
\item $G(e,2^{e-5}-1)+G(e+1,2^{e-5})+G(e+2,2^{e-4})$;
\item $G(e,2^{e-4}-1)+G(e+1,0)+G(e+2,2^{e-4})$.
\end{itemize}
\end{rmk}

Before proving the theorem, we provide the easy deduction of \fullref{1better}.
\begin{proof}[Proof of \fullref{1better}] From \fullref{gdsum} and \eqref{complement},  we obtain
\begin{equation}\sum_{i=0}^{t-1} g(e+i,k_i)\le (2^e-8)-(e-2)=2^e-e-6.\label{ge}\end{equation}
Let $e$ be fixed, and for $t\ge1$ and $0\le\ell\le 2^{e-3}-1$, let
$$M(t,\ell)=\max\biggl(\sum_{i=0}^{t-1} g(e+i,k_i):\ k_0+\cdots+k_{t-1}=\ell\biggr).$$
Then $M(t,\ell)=\max(M(t-1,i)+g(e+t-1,\ell-i):\ 0\le i\le\ell)$. Using \fullref{compatprop}, induction on $t$, and \fullref{gthm}.(6),
we obtain that for all $t$ and $\ell\le 2^{e-3}-1$
$$\gd((2^e+\cdots+2^{e+t-1}),\ell)\le M(t,\ell)$$
compatibly for all $\ell$. By \eqref{ge}, $M(t,2^{e-3}-1)\le 2^e-e-6$. Since $\gd(n,k)=\gd(n\xi_{8k+7})$, we obtain the conclusion
of \fullref{1better}.
\end{proof}

The proof of \fullref{gdsum} is expedited by the following lemma.
\begin{lem} Let $d>0$. If $G(e,i)<S(i)$ and $G(e+d,j)<S(j)$, then
$$G(e,i)+G(e+d,j)\ge G(e,i+j).$$\label{d>0}
\end{lem}
\begin{proof} This follows exactly as in the proofs of Cases 1, 4, and 5 of
\fullref{not8k+1} and the proof of \fullref{1mod8}. In those results, there was an
extra 1 on the LHS, but the larger $e$-components here more than compensate for that.
\end{proof}
\begin{rmk}{\rm Lemma \ref{d>0} is not always true when $S(-)$ is involved.
For example, if $e\ge15$, then $G(e,2)+G(e+1,3)=12<13=G(e,5)$.}
\end{rmk}
\begin{proof}[Proof of \fullref{gdsum}] Let $\cal S$ denote the set of those $k_i$
for which $G(e+i,k_i)=S(k_i)$. This includes cases in which $k_i=0$ or $k_i=1$.
If $\cal S$ is empty, then the result follows by induction from
\fullref{d>0}, since $G(e,2^{e-3}-1)=e-2$. Let $K=\sum_{k_i\in\cal S}k_i$. We split the LHS of \eqref{immsum} as
\begin{equation}\sum_{k_i\in\cal S}G(e+i,k_i)+\sum_{k_i\not\in\cal S}G(e+i,k_i).\label{splitsum}\end{equation}
Since, as is easily proved, $S(k)\ge\frac32k$, the first half of \eqref{splitsum} is $\ge\frac32K$, while \fullref{d>0}
implies that the second half of \eqref{splitsum} is $$\ge G(e,2^{e-3}-1-K)\ge e-6-\nu'(K+1)-D(e,2^{e-3}-1-K),$$
where $D(-,-)$ denotes the deviation, i.e., the sum in \eqref{geq1}. Now the desired inequality reduces to
\begin{equation}\label{32}\tfrac 32K\ge \nu'(K+1)+4+D(e,2^{e-3}-1-K).\end{equation}
If $K\ne1,3$, this inequality is true, usually with much to spare. Indeed, $K\ge\nu'(K+1)+4$ if $K\ne1,3$,
and \begin{equation}\label{34}\tfrac12K\ge D(e,2^{e-3}-1-K).\end{equation}
To see \eqref{34}, note that for $D(e,k)$ to be positive due to a single deviation, then $k=2^t\a+k_0$ with
$t> R(k_0)>k_0$, $\a$ odd, and $D(e,k)=t-R(k_0)$. For such $k$, if $k=2^{e-3}-1-K$, then $K\ge 2^t-1-k_0$, and so the difference
in \eqref{34} is $$\ge\tfrac12(2^t-1-k_0)-(t-R(k_0))=(\tfrac12(2^t-1)-t)+(R(k_0)-\tfrac12k_0)>0,$$
and a similar analysis applies when multiple deviations are involved. When $K=1,3$, \eqref{32} is true if the LHS is replaced by $S(K)=8$.
\end{proof}

\def\line{\rule{.6in}{.6pt}}

\end{document}